\documentclass[twocolumn]{autart}    
\makeatletter
\long\def\@maketablecaption#1#2{\@tablecaptionsize
    \global \@minipagefalse
    \hbox to \hsize{\parbox[t]{\hsize}{\centering #1 \\ #2}}}

\makeatother

\usepackage{amsmath,amsfonts,amssymb,amscd,latexsym,enumerate,theorem,bbm}

\usepackage{graphicx,epsfig,psfrag}
\usepackage{subfigure}
\usepackage{tikz}
\usepackage{tkz-berge}
\usepackage{cite}
\usepackage{dsfont}
 \allowdisplaybreaks
\usepackage{algorithmic}
\usepackage{float,booktabs}
\usepackage{graphicx,epsfig}
\usepackage{epstopdf}
\usepackage{siunitx}
\usepackage{tikz}
\usepackage{tkz-berge}
\usetikzlibrary{positioning}
\usetikzlibrary{arrows,%
                petri,%
                topaths}%

\usetikzlibrary{decorations.markings}         
\tikzstyle{vertex}=[circle, shading = ball, ball color = white!100!white, minimum size = 15pt, draw, inner sep=0pt]  

\newcommand{\bbm}{\begin{bmatrix}}
\newcommand{\ebm}{\end{bmatrix}}
\newcommand{\R}{\ensuremath{\mathbb R}}
\newcommand{\C}{\ensuremath{\mathbb C}}
\DeclareMathOperator{\col}{col}
\DeclareMathOperator{\im}{im}
\DeclareMathOperator{\diag}{diag}
\DeclareMathOperator{\bdiag}{blockdiag}
\usepackage{subeqnarray}

\usepackage{epsfig}

\newcommand{\calA}{\ensuremath{\mathcal{A}}}
\newcommand{\calB}{\ensuremath{\mathcal{B}}}
\newcommand{\calC}{\ensuremath{\mathcal{C}}}

\newcommand{\calP}{\ensuremath{\mathcal{P}}}

\newcommand{\calU}{\ensuremath{\mathcal{U}}}

\newcommand{\calG}{\ensuremath{\mathcal{G}}}
\newcommand{\calV}{\ensuremath{\mathcal{V}}}
\newcommand{\calE}{\ensuremath{\mathcal{E}}}
\newcommand{\calN}{\ensuremath{\mathcal{N}}}

\newcommand{\bse}{\begin{subequations}}
\newcommand{\ese}{\end{subequations}}
\def\be{\begin{equation}}
\def\ee{\end{equation}}
\newcommand{\EW}{\hspace*{\fill} $\square$\noindent}

\usepackage[prependcaption,colorinlistoftodos]{todonotes}

{\theorembodyfont{\itshape}\newtheorem{theorem}{Theorem}}
{\theorembodyfont{\itshape}\newtheorem{proposition}[theorem]{Proposition}}
{\theorembodyfont{\itshape}\newtheorem{lemma}[theorem]{Lemma}}
{\theorembodyfont{\itshape}}
{\theorembodyfont{\upshape}\newtheorem{definition}[theorem]{Definition}}
{\theorembodyfont{\itshape}}
{\theorembodyfont{\upshape}\newtheorem{remark}[theorem]{Remark}}
{\theorembodyfont{\upshape}}
{\theorembodyfont{\upshape}\newtheorem{example}[theorem]{Example}}
{\theorembodyfont{\upshape}\newtheorem{assumption}{Assumption}}

\def\QEDopen{{\setlength{\fboxsep}{0pt}\setlength{\fboxrule}{0.2pt}\fbox{\rule[0pt]{0pt}{1.3ex}\rule[0pt]{1.3ex}{0pt}}}}
\def\QED{\QEDopen} 

\newenvironment{pfof}[1]{\vspace{1ex}\noindent{\itshape Proof of
    #1:}\hspace{0.5em}} {\hfill\QED\vspace{1ex}}

\newcommand{\BP}{\noindent{\bf Proof. }}
\newcommand{\EP}{\hspace*{\fill} $\blacksquare$\smallskip\noindent}

\newcommand{\ones}{\mathds{1}}

\newcommand{\part}[2]{\ensuremath{\frac{\partial #1}{\partial #2}}}

\newcommand{\norm}[1]{\ensuremath{\left\| #1 \right\|}}

\usepackage{etoolbox}
\preto\subequations{\ifhmode\unskip\fi}

\begin{document}
		
\begin{frontmatter}

\title{Secant and Popov-like Conditions in Power Network Stability}


\author[RUG]{Nima Monshizadeh}\ead{n.monshizadeh@rug.nl},$\,$    
\author[CAM]{Ioannis Lestas}\ead{icl20@cam.ac.uk}          
\address[RUG]{Engineering and Technology Institute Groningen, University of Groningen, Nijenborgh 4, Groningen, 9747 AG , The Netherlands}
\address[CAM]{Department of Engineering, University of Cambridge, Trumpington Street, Cambridge, CB2 1PZ, United Kingdom}
\thanks{This work was supported by ERC starting grant 679774.}
\begin{keyword}                           
Power network stability, Secant conditions, Popov criterion, Passivity
\end{keyword}                             

\begin{abstract}
The problem of decentralized frequency control in power networks has received an increasing attention in recent years due to its significance in modern power systems and smart grids. Nevertheless, generation dynamics including turbine-governor dynamics, in conjunction with nonlinearities associated with generation and power flow, increase significantly the complexity in the analysis, and are not adequately addressed in the literature. In this paper we show how incremental secant gain conditions can be used in this context to deduce decentralized stability conditions with reduced conservatism.
Furthermore, for linear generation dynamics, we establish Popov-like conditions that are able to reduce the conservatism even further by incorporating additional local information associated with the coupling strength among the bus dynamics. Various examples are discussed throughout the paper to demonstrate the significance of the results presented.
\end{abstract}
\end{frontmatter}
\section{Introduction}

{
Due to the large scale penetration of renewable energy sources in the power grid, there has been an increasing interest in recent years in decentralized and distributed frequency control schemes in power networks. 
As a result of the nonlinearities associated with power flows, and also potential nonlinearites in the generation dynamics, a Lyapunov analysis is a natural tool often used in this context for stability analysis, e.g. the use of energy functions in \cite{tsolas1985structure,chu1999constructing,chang1995direct},  and more recent Lyapunov approaches in \cite{de2018bregman,trip2016internal,zhao2015distributed,jayawardhana-passivity,arcak2016networks}.
Extensions to differential algebraic models can be found in
\cite{schiffer2016stability,de2016lyapunov}, see also \cite{Hill-DAE}.

Nevertheless, a feature that can complicate significantly such a Lyapunov analysis is the presence of turbine/governor dynamics in conjunction with nonlinearities often present in the generation or controllable demand side, such as deadbands and saturation. Such dynamics are often not explicitly addressed in the literature and various notable exceptions either resort to linearizations or propose gain conditions relative to the system damping that are lower than those encountered in practical implementations \cite{trip2016optimal,zhao2014optimal}. In \cite{kasis2017primary} a passivity property on the aggregate bus dynamics was proposed, with further generalizations provided in \cite{devane2017primary}, as a means of reducing the conservatism in the analysis. Systematic methods exist for verifying these properties for the case of linear systems. However, in the case of nonlinear systems, the problem of determining the minimum damping needed to passivate it is in general a nontrivial problem, as the form of the underlying storage function is unknown. Furthermore, simpler approaches that achieve passivation via a restricted $\mathcal{L}_2$ gain, can in general be restrictive, hence alternative methodologies need to be investigated.

%
In this paper, we show that the use of suitable incremental secant conditions, inspired by
\cite{sontag2006passivity,arcak2006diagonal}, can facilitate the construction of classes of Lyapunov functions in this context and lead to stability conditions with reduced conservatism, in cases where a linearizaton is not appropriate. These conditions are decentralized and result in asymptotic stability for a range of equilibria of the system, and thus are able to cope with the uncertainties in load parameters and generation setpoints. The applicability of these conditions is demonstrated by means of several examples, which illustrate that these provide stability guarantees with larger control gains, which in turn will enhance the performance of the network.

Furthermore, for the case where generation dynamics are linear,
we show how even less conservative stability conditions can be established by
leveraging additional local information associated with the system model.
In particular,
we maintain the nonlinearity of the power flows in the analysis, and use Popov like arguments to derive distributed conditions
that take into account the coupling strength among the bus dynamics.
Numerical examples are also used to investigate the relative merits of the conditions derived. 
It is observed that while the secant conditions, and more generally conditions relying on passive bus dynamics, are suitable for strongly coupled networks, the latter Popov-like conditions can offer 
improvement when the bus dynamics are more weakly coupled.

The structure of the paper is as follows. The power network model is provided in Section \ref{s:dae}. The desired asymptotic behavior of the system is characterized in Section \ref{s:synch}. The main results of the paper are provided in Section \ref{s:main}, and several examples are discussed to illustrate the applicability of the proposed stability conditions. The paper closes with conclusions in Section \ref{s:conclusions}.

\noindent{\textbf{Notation.}} The $n \times n$ identity matrix is denoted by $I_n$, and $\ones_n$ is the vector of all ones in $\R^n$, where the subscript is dropped if no confusion may arise. For $i\in \{1,2, \ldots, n\}$, by $\col(a_i)$ we denote the  vector $(a_1, a_2, \ldots, a_n)$. For given vectors $a\in \R^n$ and $b\in \R^m$, we denote the vector $( a^{\sf T}, b^{\sf T})^{\sf T}\in \R^{n+m}$ by $\col(a, b)$ or sometimes simply by $(a, b)$. Given  a map ${H}:\R^n \to \R$, its transposed gradient is denoted by
$\nabla{H}:= \left(\frac{\partial {H}}{\partial x}\right)^{\sf T}$.
\section{Differential-Algebraic model of power network}\label{s:dae}

We consider a structure-preserving model of power networks composed of load and generation buses.
The topology of the grid is represented by a connected and undirected graph $\calG(\calV,\calE)$ with a vertex set (or buses) $\calV=\{0,1, \ldots, n\}$, and an edge set $\calE$ given by the set of unordered pairs $\{i, j\}$ of distinct vertices $i$ and $j$.  The cardinality of $\calE$ is denoted by $m$.
%
We assume that the line admittances are purely inductive, and two nodes $\{i, j\}\in \calE$ are connected by a nonzero real susceptance $\beta_{ij}<0$.
The set of neighbors of the $i^{th}$ node is denoted by $\calN_i=\{j\in \calV \mid {\{i,j\}\in \calE}\}.$ The voltage phase angle at node $i\in \calV$ is denoted by $\theta_i\in \R$. Voltage magnitudes $V_i\in \R^+$ are assumed to be constant.

The set of generators is given by $\calV_g=\{0,1, \cdots, n_g\}$. For each generator $i\in \calV_g$,
the phase $\theta_i$ evolves according to  \cite{kundur94}
\bse\label{e:theta}
\begin{align}
\dot\theta_i&=\omega_i\\
\label{e:swing}
M_i \dot{\omega}_i&=-D_i\omega_i-p_{i}(\theta)+p_i^\ast+u_i\,,
\end{align}
\ese
where
\be\label{e:active}
p_i(\theta)=\sum_{j\in \calN_i}\nolimits |\beta_{ij}|V_iV_j \sin(\theta_i-\theta_j)
\ee
is the active power drawn from bus $i$. Here, $\omega_i$ is the frequency deviation from the nominal frequency (namely \SI{50}{\Hz}), $M_i>0$ is the inertia constant, $D_i>0$ is the damping constant, the constant $p_i^\ast$ is the active power setpoint, and $u_i\in \R$ is the additional local power generation at bus $i$.  The constant $p_i^*$ may also capture the constant power loads collocated with the $i$th generator bus.

As for the loads, we consider constant power loads given by algebraic equations
\be\label{e:loads}
0=p_i^*-p_i(\theta)\,,
\ee
for each $i\in \calV_\ell=\calV\setminus \calV_g$, where $p_i(\theta)$ is given by \eqref{e:active} and $p_i^*$ is constant. Note that constant impedance loads behave similarly to constant power loads if the voltages are approximately constant.
We remark that the exact value of $p_i^*$, $i\in \calV$, is not known a priori.

To capture a broad class of generation dynamics, let $u_i\in \R$ be given by a nonlinear system of the form
\bse\label{e:ui}
\begin{align}
\dot \xi_i&= f_i(\xi_i, -\omega_i)\\
u_i&=h_i (\xi_i,-\omega_i)\,,
\end{align}
\ese
where $f_i: \R^{n_i} \times \R \rightarrow \R^{n_i}$ and $h_i: \R^{n_i}\times \R\rightarrow  \R$ are continuous and locally Lipschitz.
We sometimes denote such dynamical systems by $\Sigma_i(-\omega_i,\xi_i,u_i)$ in short.
For any constant input $\omega_i=\overline \omega_i$, we assume that \eqref{e:ui} possesses an isolated equilibrium $\xi_i=\overline \xi_i$, and we write $\overline u_i=h_i(\overline \xi_i, -\overline \omega_i)$. We also assume that such an equilibrium is observable from the constant input-output pair $(-\overline \omega_i, \overline u_i)$, i.e.,
$\dot \xi_i= f_i(\xi_i, -\overline\omega_i)$ together with $h_i(\xi_i, -\overline \omega_i)=h_i(\overline \xi_i, -\overline \omega_i)$ implies that $\xi_i=\overline \xi_i$.
Note that the dynamics \eqref{e:ui} may include primary control, controllable loads, turbine governor dynamics, and possible static nonlinearities  in the generation dynamics. Examples of higher order turbine-governor dynamics include models for steam turbines (with or without reheat) as in e.g.  \cite[Sec. 11.1.4]{kundur94}, \cite[Sec. 11.3.1]{Machowski2008}.

The power network dynamics can be written in vector form as the following differential algebraic system:
\bse\label{e:theta-compact}
\begin{align}
\dot\theta_g&=\omega_g \\
M \dot{\omega}_g&=-D\omega_g-p_g(\theta)+p_g^\ast+h(\xi, -\omega_g)\\
\dot \xi&= f(\xi, -\omega_g)\\
0&=-p_\ell(\theta)+ p_\ell^*\;,
\end{align}
\ese
where $M=\bdiag(M_i)$, $D=\bdiag(D_i)$, $\theta_g=\col(\theta_i)$, $\omega_g=\col(\omega_i)$, $p_g(\theta)=\col(p_i(\theta))$, $p_g^\ast=\col(p_i^*)$, $\xi=\col(\xi_i)$, $f=\col(f_i)$, and $h=\col(h_i)$ for $i\in \calV_g$. Similarly,
$p_\ell(\theta)=\col(p_i(\theta))$ and $p_\ell^\ast=\col(p_i^*)$, $i\in \calV_\ell$.

Let $R$ be the incidence matrix of the graph. Note that, by associating an arbitrary orientation to the edges, the incidence
matrix $R\in\R^{(n+1) \times m}$ is defined element-wise as $R_{ik} = 1,$ if node
$i$ is the sink of the edge $k$, $R_{ik} = -1,$ if $i$ is the source of the edge $k$, and $R_{ik} = 0$ otherwise.
In addition, let $\Gamma:=\diag(\gamma_k), \gamma_k=|\beta_{ij}|V_iV_j,$ for each edge $k \sim \{i, j\}$ of the graph, where the edge numbering is in agreement with the incidence matrix $R$.
Then the vector of active power transfer $p(\theta)=\col(p_g(\theta), p_\ell(\theta))$ is written as
\be\label{e:active-vector}
p(\theta)=R\Gamma \boldsymbol\sin(R^{\sf T}\theta)=\bbm R_g \\ R_\ell\ebm \Gamma \boldsymbol\sin(R^{\sf T}\theta)\,,
\ee
where $R_g$ and $R_\ell$ are the  {submatrices} of $R$ obtained by collecting the rows of $R$ indexed by $\calV_g$ and $\calV_\ell$, respectively.
The operator $\boldsymbol\sin(\cdot)$ is interpreted element-wise.

\section{Synchronous solution and a change of coordinates}\label{s:synch}

We are interested in a {\em synchronous motion } of the power network, where the voltage phasors rotate with the same frequency.
This writes as $\overline\theta_i(t)=\omega^*t+\overline \theta_i(0)$ for each $i\in \calV$, with  constant $\overline\theta_i(0)\in \R$.  Note that a synchronous motion explicitly depends on time.
In addition, note that if $(\theta, \omega_g, \xi)$ is a solution to \eqref{e:theta-compact}, then $(\theta+ c\ones_{n+1}, \omega_g, \xi)$ is also a solution to \eqref{e:theta-compact} for any constant $c\in \R$. To get around this rotational invariance, we perform a change of coordinates by taking a phase angle of a generation bus, namely $\theta_0$, as a reference:
 \be\label{varphi}
\varphi_i = \theta_i -\theta_0\,, \quad i=1,\ldots, n.
\ee
This new set of coordinates satisfies
\[
\begin{bmatrix}
0\\
\varphi_1\\
 \vdots \\ \varphi_{n}
\end{bmatrix}=
\begin{bmatrix}
\theta_0\\
\theta_1 \\
\vdots \\ \theta_n
\end{bmatrix}-\mathds{1}_{n+1} \theta_0.
\]
Let $R_\varphi\in \R^{n\times m}$ denote the incidence matrix with its first row removed, and let $\col(\varphi_i):=\varphi\in \R^n$. Then, by the equality above and noting that $\ones \in \ker R^{\sf T}$, we have
$$R^{\sf T}\theta=R_\varphi^{\sf T} \varphi.$$
Moreover, we have $\varphi= E^{\sf T} \theta$ where $E^{\sf T}=\bbm -\ones_n & I_n\ebm$.
This can be rewritten as $\varphi=E_g^{\sf T}\theta_g+E_\ell^{\sf T}\theta_\ell$, where the matrix $E$ is partitioned accordingly as $E^{\sf T}=\bbm E_g^{\sf T} & E_\ell^{\sf T}\ebm$.
Now, let $\varphi_g:=E_g^{\sf T} \theta_g$ and $\varphi_\ell:= E_\ell^{\sf T}\theta_\ell$.
To clarify note that $\varphi, \varphi_g, \varphi_\ell \in \R^n$ and $\varphi=\varphi_g+\varphi_\ell$.

Then, the system \eqref{e:theta-compact} in the new coordinates reads as
\bse\label{e:varphi-compact0}
\begin{align}
\dot\varphi_g&=E_g^{\sf T} \omega_g\\
M \dot{\omega}_g&=-D\omega_g-R_g\Gamma \boldsymbol{\sin}(R_\varphi^{\sf T}\varphi)+p_g^\ast+h(\xi, -\omega_g)\\
\dot \xi&= f(\xi, -\omega_g)\\
0&=-R_\ell\Gamma \boldsymbol{\sin}(R_\varphi^{\sf T}\varphi)+ p_\ell^*\,.
\end{align}
\ese

Let $U(\varphi):=-\ones_m^{\sf T}\Gamma \boldsymbol{\cos}(R_\varphi^{\sf T}\varphi)$, where again $\boldsymbol{\cos}(\cdot)$ is defined element-wise.
Clearly, $\nabla U(\varphi)= R_\varphi \Gamma \boldsymbol{\sin}(R_\varphi^{\sf T}\varphi)$.  In addition, it is easy to see that $R=ER_\varphi$, and thus $R_g=E_gR_\varphi$, $R_\ell=E_\ell R_\varphi$.
{Then, \eqref{e:varphi-compact0}} can be written as
\bse\label{e:varphi-compact}
\begin{align}
\dot\varphi_g&=E_g^{\sf T} \omega_g\\
M \dot{\omega}_g&=-D\omega_g- E_g \nabla U(\varphi)+p_g^\ast+h(\xi, -\omega_g)\\
\dot \xi&= f(\xi, -\omega_g)\\
0&=-E_\ell\nabla U(\varphi)+ p_\ell^*\,.
\end{align}
\ese


The representation above gives a differential algebraic model of the form
\bse\label{e:dae-generic}
\begin{align}
\dot x&= F(x, q)\\
0&=g(x, q)\,,
\end{align}
\ese
where $x=\col(\varphi_g, \omega_g, \xi)$ and $q=\varphi_\ell$, noting that $\varphi=\varphi_g+\varphi_\ell$.
We assume that initial conditions are compatible with the algebraic equations, i.e., $0=g(x(0), q(0))$. For now, we also assume that the system above has a unique solution, for a nonzero interval of time, starting from any compatible initial condition. As will be observed later, this assumption is automatically satisfied since we will work in a region of state space where the algebraic constraints are {\em regular}, i.e., $\frac{\partial g}{\partial q}$ has full row rank.

As a result of this change of coordinates, a synchronous motion of the power network will be mapped to an equilibrium of the differential algebraic system \eqref{e:varphi-compact}, namely the point $(\overline \varphi, \overline \omega_g, \overline \xi)$ where $\overline \omega_g=\ones \omega^*$ with $\omega^*\in \R$ being constant, and $\overline \varphi\in \R^n$ and  $\overline \xi \in \R^N$, with $N=\sum_{i\in \calV_g} n_i$, are constant vectors satisfying
\bse\label{e:feas}
\begin{align}
0&=-D\ones \omega^* - E_g \nabla U(\overline\varphi)+p_g^\ast+h(\overline\xi, -\ones \omega^*),\\
0&=-E_\ell\nabla U(\overline \varphi)+ p_\ell^*\\
0&=f(\overline \xi, -\ones \omega^*)\,.
\end{align}
\ese

We will refer to the equilibrium point $(\overline \varphi, \overline \omega_g, \overline \xi)$ as the {\em synchronous solution} of the power network.
By \eqref{e:feas}, existence of such a solution imposes the following feasibility assumption:
\begin{assumption}\textbf{(Existence of a synchronous solution)}\label{a:feas}
There exists a constant $\omega^*\in \R$, a constant vector $\overline \varphi\in \R^{n}$ with $R_\varphi^{\sf T}\overline\varphi \in (-\frac{\pi}{2}, \frac{\pi}{2})^m,$ and a constant vector $\overline \xi\in \R^{N}$ such that \eqref{e:feas} is satisfied.
\end{assumption}

The additional requirement that $R_\varphi^{\sf T}\overline \varphi \in (-\frac{\pi}{2}, \frac{\pi}{2})^m$ means that the relative phase angles at steady-state should belong to the interval $(-\frac{\pi}{2}, \frac{\pi}{2})$. This assumption is often referred to as the security constraint and is ubiquitous in the literature, see e.g. \cite{trip2016internal,dorfler2016breaking,kasis2017primary}. In case of linear generation dynamics, the description of \eqref{e:feas} and consequently Assumption \ref{a:feas} can be made more explicit, see Lemma \ref{l:equib-lin} and Assumption \ref{a:feas-lin}.
\section{Main results}\label{s:main}

\subsection{Incremental passivity of the differential algebraic model}
Consider  the differential algebraic system
\bse\label{e:dae-open}
\begin{align}
\dot\varphi_g&=E_g^{\sf T} \omega_g\\
M \dot{\omega}_g&=-D\omega_g- E_g \nabla U(\varphi)+p_g^\ast+u\\
0&=-E_\ell\nabla U(\varphi)+ p_\ell^*\\
y&=\omega_g
\end{align}
\ese
%
with input-state-output $(u, (\varphi, \omega_g), \omega_g)$, where $u=\col(u_i)$.
Clearly, \eqref{e:varphi-compact} can be seen as a negative feedback interconnection of  \eqref{e:dae-open} with \eqref{e:ui}.
As a first step towards a systematic stability analysis of \eqref{e:theta-compact}, we identify an incremental passivity property of \eqref{e:dae-open}
with respect to a synchronous solution $(\overline \varphi, \overline \omega_g, \overline \xi)$. To formalise this property, we need the following definition:

\begin{definition}\label{d:passive}
Consider the differential algebraic system
\bse\label{e:dae}
\begin{align}
\dot x_o&=f(x_o, x_a, u)\\
0&=g(x_o,x_a)\\
y&=h(x_o, u)
\end{align}
\ese
with input-state-output $(u, x, y)$, where $x=\col(x_o, x_a)$.
System \eqref{e:dae}
is incrementally passive  with respect to a point $(\overline u, \overline x ,\overline y)\in \mathcal{U} \times \mathcal{X} \times \mathcal{Y}$, with $\overline y=h(\overline x_0, \overline u)$, if there exists a nonnegative\footnote{Nonnegativity is assumed in $\mathcal{X}$. The set $\mathcal{X}$ can be shrunk as desired.} and continuously differentiable function $S(x)$
and a positive semidefinite matrix $Q$, such that
for all $x\in \mathcal{X}$ and $u\in \calU$,
the inequality
\be\label{ine:passive}
\dot S(x)\le
-(y-\overline y)^{\sf T} Q(y-\overline y)+(y-\overline y)^{\sf T} (u-\overline u)
\ee
holds. In case the matrix $Q$ is positive definite, we call the system output strictly incrementally passive with respect to
$(\overline u, \overline x, \overline y)$.
\end{definition}

\medskip{}
Now, we have the following proposition:

\bigskip{}
\begin{proposition}\label{p:passivity}
Let $(\overline \varphi, \overline \omega_g)$, with $R_\varphi^{\sf T} \overline \varphi \in (-\frac{\pi}{2}, \frac{\pi}{2})^n$, be an equilibrium of \eqref{e:dae-open} for some constant input $u=\overline u$, and let $\overline y=\overline \omega_g$. Then the differential algebraic system \eqref{e:dae-open} is output strictly incrementally passive with respect to $(\overline u, (\overline\varphi, \overline \omega_g), \overline y)$. {In particular, the storage function $S(\varphi, \omega_g)$ given by \eqref{e:Storage} satisfies}
\be\label{e:passivity-ineq}
\dot S = -(\omega_g- \overline \omega_g)^{\sf T} D (\omega_g -\overline \omega_g) +  (\omega_g- \overline \omega_g)^{\sf T}(u-\overline u)\,.
\ee
{Moreover, this storage function has a local strict minimum at $(\overline \varphi, \overline \omega_g)$.}
\end{proposition}
\BP
Noting that $(\overline \varphi, \overline \omega_g)$ is an equilibrium of \eqref{e:dae-open}, we can rewrite \eqref{e:dae-open} as
\bse\label{e:dae-open-inc}
\begin{align}
\dot \varphi_g&=E_g^{\sf T} \omega_g\\
\label{dae-open-inc-wg}
\nonumber
M \dot{\omega}_g&=-D(\omega_g-\overline \omega_g)\\ & \qquad  - E_g \big(\nabla U(\varphi)-\nabla U(\overline\varphi)\big)+ u-\overline u\\
\label{dae-open-inc-alg}
0&=-E_\ell \big(\nabla U(\varphi)-\nabla U(\overline \varphi)\big)\\
y&=\omega_g\,.
\end{align}
\ese
Take the storage function candidate
\begin{align}\label{e:Storage}
\nonumber
&S=\frac{1}{2}(\omega_g-\overline \omega_g)^{\sf T}M(\omega_g-\overline \omega_g)\\
&+   U(\varphi) -  U(\overline \varphi)-(\varphi - \overline \varphi)^{\sf T}  \nabla U(\overline \phi)\;. 
\end{align}
The first term of $S$ is clearly nonnegative and is equal to zero whenever $\omega_g=\overline \omega_g$.
The terms in the second line of \eqref{e:Storage} constitute a Bregman distance defined for the function $U(\varphi)$ with respect to the point $\varphi=\overline \varphi$, \cite{bregman1967relaxation,jayawardhana-passivity,de2018bregman}.
Note that $\nabla^2U(\varphi)= R_\varphi \Gamma [\boldsymbol{\cos}(R_\varphi^{\sf T} \varphi)] R_\varphi^{\sf T}$, and that $R_\varphi$ has full row rank. Here,  $[\boldsymbol{\cos}(R_\varphi^{\sf T} \varphi)]$ denotes the diagonal matrix constructed from the vector $\boldsymbol\cos(R_\varphi^{\sf T} \varphi)$.
Hence, we find that $U$ is a strict convex function of $\varphi$, as long as the relative phase angles $R_\varphi^{\sf T} \varphi$ belong to a closed subset $\mathcal{X}_\varphi$ of
$(-\frac{\pi}{2},\frac{\pi}{2})^m$. Consequently, the aforementioned Bregman distance is strictly positive whenever $\varphi \neq \overline \varphi$ and $R_\varphi^{\sf T} \varphi \in \mathcal{X}_\varphi$.
Moreover, the partial derivatives of $S$ are computed as
\be\label{e:partial}
\frac{\partial S}{\partial\omega_g}=\omega_g-\overline \omega_g, \quad \frac{\partial S}{\partial\varphi}=\nabla U(\varphi)- \nabla U(\overline \varphi).
\ee
Therefore the partial derivatives of $S(\varphi, \omega_g)$ vanish at $(\overline \varphi, \overline \omega_g)$, and thus $S$ has a local strict minimum at this point.

Noting that $\varphi=\varphi_g+\varphi_\ell$, the partial derivative of $E_L\nabla U(\varphi)$ with respect to the state variable associated to the algebraic equations, namely $\varphi_g$, is obtained as $E_\ell \nabla ^2 U(\varphi)$. This matrix has full row rank since the matrix $E_\ell$ has full row rank and $U(\varphi)$ is strictly convex in the region for which $R_\varphi^{\sf T} \varphi \in \mathcal{X}_\varphi$. Therefore, starting from a compatible initial condition,
there exists a unique solution $((\varphi_g, \varphi_\ell), \omega_g)$ satisfying the differential algebraic equations \eqref{e:dae-open-inc}, for some nonzero interval of time \cite{Hill-DAE}. Taking the time derivative of $S$ along such a solution yields
\begin{align*}
\dot S&= (\omega_g-\overline \omega_g)^{\sf T} \big( - D(\omega_g-\overline \omega_g) - E_g(\nabla U(\varphi)- \nabla U(\overline \varphi))\big)\\
 &\quad\;  +  (\omega_g-\overline \omega_g)^{\sf T}(u-\overline u)+(\nabla U(\varphi)- \nabla U(\overline \varphi))^{\sf T}E_g^{\sf T} \omega_g\,.
\end{align*}
This simplifies to
\begin{align}\label{e:dotS-proof}
\nonumber
\dot S&= -(\omega_g-\overline \omega_g)^{\sf T}D(\omega_g-\overline \omega_g)+   (\omega_g-\overline \omega_g)^{\sf T}(u-\overline u)\\
& \qquad \quad +\overline \omega_g^{\sf T} E_g (\nabla U(\varphi)- \nabla U(\overline \varphi))\,.
\end{align}
Note that $\overline \omega_g \in \im \ones_{n_g}$, namely $\overline  \omega_g=\ones_{n_g} \omega^*$ for some constant $\omega^*\in \R$.
Then, it is easy to see that
\[
\overline \omega_g^{\sf T} E_g (\nabla U(\varphi)- \nabla U(\overline \varphi))= -\omega^* \ones^{\sf T} E_\ell (\nabla U(\varphi)- \nabla U(\overline \varphi))\,.
\]
The right hand side of the equality above is equal to zero by the algebraic equation \eqref{dae-open-inc-alg}, and therefore \eqref{e:dotS-proof} reduces to \eqref{e:passivity-ineq}.
\EP

\bigskip{}
Recall that \eqref{e:varphi-compact} is given by a negative feedback interconnection of  \eqref{e:dae-open} with \eqref{e:ui}.
In case the generation dynamics \eqref{e:ui} is (incrementally) passive as well, then by exploiting the result of Proposition \ref{p:passivity}, the closed-loop system enjoys suitable stability properties due to the standard results on interconnection of passive systems, see Example \ref{e:simple}. On the other hand, if \eqref{e:ui} is not (incrementally) passive, then stability of closed-loop system is not automatically guaranteed, and it requires additional conditions.
In particular, the droop/control gain {needs to be restricted, see e.g.} \cite[Ex. 11.3]{bergen1999power}.
\medskip{}
\begin{example}\label{e:simple}
Suppose that the generation dynamics are given by static input-output relation $u_i=h_i(-\omega_i)$, where $h_i$ is a strictly increasing map for each $i\in \calV_g$. Then, clearly,
$(\omega_g- \overline \omega_g)^{\sf T} (u- \overline u) \leq 0$.  Substituting this into \eqref{e:passivity-ineq} concludes stability of the equilibrium
$(\overline \varphi, \overline \omega_g, \overline \xi)$, {as $S$ has a strict minimum at $(\overline \varphi, \overline \omega_g)$ and} $\dot S$ is nonpositive. Asymptotic stability follows by a suitable analysis of an invariant set of the system. 
\EW
\end{example}
\subsection{Small incremental-gain conditions}

Considering the nonlinearity of the generation dynamics, a first approach is to use an incremental $L_2$-gain {{argument}.
First, the following definition is needed:
\begin{definition}({\textbf{Incremental $L_2$ stablility)}}
The system
\bse
\begin{align}
\dot x&=f(x, u)\\
y&=h(x, u)
\end{align}
\ese
with input-state-output $(u, x, y)$ is incrementally $L_2$ stable with respect to a point $(\overline u, \overline x ,\overline y)\in \mathcal{U} \times \mathcal{X} \times \mathcal{Y}$, with $\overline y=h(\overline x, \overline u)$, if there exists a nonnegative continuously differentiable function $S(x)$
and a scalar $\delta\in \R^+$  such that
for all $x\in \mathcal{X}$ and $u\in \calU$,
the inequality
\be\label{ine:L2}
\dot S(x)\le  -\norm{y-\overline y}^2+ \delta^2 \norm{u-\overline u}^2
\ee
holds. The system has an incremental $L_2$-gain not greater than $\delta$ in this case.
\end{definition}
\bigskip{}

The notions of stability and asymptotic stability used here are those of \cite{Hill-DAE}. Now, we have the following small incremental-gain result:
\bigskip{}

\begin{proposition}\label{p:L2}
Let Assumption \ref{a:feas} hold. Assume that \eqref{e:ui} is incrementally $L_2$ stable with respect to $(-\overline\omega_i, \overline \xi_i, \overline u_i)$ and {that the associated storage function has a strict minimum at this point}. Let the incremental $L_2$-gain of \eqref{e:ui} be not greater than
$\delta_i$ for each $i\in \calV_g$.
Then,  $(\overline \varphi, \overline \omega_g, \overline \xi)$ is an asymptotically stable equilibrium of \eqref{e:varphi-compact} if, for each $i$,
\be\label{e:L2-condition}
\delta_i <D_i\,.
\ee
\end{proposition}
\medskip{}

\BP
By \eqref{e:passivity-ineq}, it is easy to verify that
\begin{align*}
&2\dot S =\\
&\quad  - \big(u-\overline u - D(\omega_g- \overline \omega_g)\big)^{\sf T} D^{-1} \big(u-\overline u - D(\omega_g- \overline \omega_g)\big)\\
&\quad - (\omega_g-\overline \omega_g)^{\sf T} D (\omega_g-\overline \omega_g) + (u-\overline u)^{\sf T} D^{-1}  (u-\overline u)\,.
\end{align*}
 Moreover, by assumption, for each $i\in\calV_g$, there exists a storage function $Z_i(\xi_i)$ with its minimum at $\xi_i=\overline \xi_i$, satisfying
\[
\dot Z_i\leq -(u_i-\overline u_i)^2 + \delta_i^2(\omega_i-\overline \omega_i)^2\,.
\]
By \eqref{e:L2-condition}, there exists $\lambda\in \R^+$ such that $\delta_i< \lambda_i < D_i$.
Now, let $$Z(\xi):=\frac{1}{2}\sum_i \frac{D_i}{\lambda_i^2}\, Z_i(\xi_i)\,.$$ Then we have
\[
2 \dot Z\leq - \sum_i \frac{D_i}{\lambda_i^2} (u_i-\overline u_i)^2  + \sum_i  \frac{D_i\delta_i^2}{\lambda_i^2} (\omega_i -\overline \omega_i)^2\,.
\]
Hence, we find that
\begin{align*}
2 \dot S+ 2\dot Z \leq  \;&\sum_i (\frac{1}{D_i}-\frac{D_i}{\lambda_i^2}) (u_i-\overline u_i)^2  \\* &+ \quad \sum_i  (\frac{D_i\delta_i^2}{\lambda_i^2}-D_i) (\omega_i -\overline \omega_i)^2.
\end{align*}
Due to the fact that $\delta_i< \lambda_i < D_i$, the right hand side of the above inequality is nonpositive, and is equal to zero if and only if
$(\omega_g, u)=(\overline \omega_g, \overline u)$. Note that $S(\varphi, \omega_g)+Z(\xi)$ has a local strict minimum at $(\overline \varphi, \overline \omega_g, \overline \xi)$. Also recall that the algebraic equations are regular in a neighborhood of this point. Then, one can construct compact level sets around  $(\overline \varphi, \overline \omega_g, \overline \xi)$ which are forward invariant. LaSalle's invariance principle with $V=S+Z$ as the Lyapunov function can then be invoked, and on any invariant set with $\dot V =0$ we have
$\omega_g=\overline \omega_g$ and $h(\xi, \overline \omega_g)=h(\overline \xi, \overline \omega_g)$. By the observability assumption of \eqref{e:ui}, we find that $\xi=\overline \xi$ on the invariant set.
Substituting this into the dynamics \eqref{dae-open-inc-wg} and \eqref{dae-open-inc-alg} yields
\begin{align*}
0&=- E_g \big(\nabla U(\varphi)-\nabla U(\overline\varphi)\big)\\
0&=-E_\ell \big(\nabla U(\varphi)-\nabla U(\overline \varphi)\big),
\end{align*}
on the invariant set. Hence, $0=E \big(\nabla U(\varphi)-\nabla U(\overline\varphi)\big)$, which noting that $E$ has full column rank results in
\[
0=\nabla U(\varphi)-\nabla U(\overline\varphi)\,.
\]
By \eqref{e:partial} and the fact that $Z$ has a strict minimum at $\overline \xi$, we observe that the partial derivatives of $S+Z$ vanish on the invariant set. Consequently, the invariant set comprises only the equilibrium
$(\overline \varphi, \overline \omega_g, \overline \xi)$, baring in mind that this point is a local strict minimum of $S+Z$. This completes the proof.
\EP
\medskip{}

\begin{example}\label{ex:2-order}
For each bus $i\in \calV_g$, let $u_i$ be given by the nonlinear second-order dynamics
\bse\label{e:example-2nd}
\begin{align}
\tau_{\alpha, i}\dot\alpha_i&=-\nabla c_i(\alpha_i) + k_i(-\omega_i)\\
\tau_{\beta, i}\dot\beta_i&= - \beta_i+\alpha_i\\
u_i&=\beta_i\,,
\end{align}
\ese
where $\alpha_i, \beta_i\in \R$ are state variables, $\tau_{\alpha, i}, \tau_{\beta, i}\in \R^+$ are time constants, and $c_i:\Omega_c \rightarrow \R$ is a strongly convex function, i.e., there exists $\rho^c_i\in \R^+$ such that
\[
(\alpha_i-\overline \alpha_i)(\nabla c_i(\alpha_i)-\nabla c_i(\overline \alpha_i))\geq \rho_i^c(\alpha_i- \overline \alpha_i)^2,
\]
for all $\alpha_i, \overline \alpha_i \in \Omega_c$. In addition, the map $k_i:\Omega_k \rightarrow \R$ satisfies $|k_i(-\omega_i)-k_i(-\overline \omega_i)|\leq  \rho^k_i |\omega_i-\overline \omega_i|$, $\forall \omega_i, \overline \omega_i\in \Omega_k$.
{For $c_i(\alpha_i)=\frac{1}{2}{\alpha_i^2}$, the  model \eqref{e:example-2nd} can represent second-order turbine governor dynamics, see e.g. \cite[Sec. 11.1]{bergen1999power}, where $k_i$ is allowed to be nonlinear in order to capture saturation, deadband, or simply a nonlinear droop gain. This case can also represent a decentralized leaky-integral controller \cite{ainsworth2013design,heidari2017ultimate,weitenberg2017robust} cascaded {with  first-order} turbine governor dynamics. Allowing for strongly convex functions other than the quadratic ones for $c_i(\cdot)$ provides additional flexibility in the design, such as more sophisticated power sharing properties compared to the proportional ones in \cite{weitenberg2017robust}.}

 Let $(\overline \alpha_i, \overline \beta_i)$ denote the equilibrium of \eqref{e:example-2nd} resulting from a constant input $-\overline \omega_i$. 
By defining $v_i=k_i(-\omega_i)$, $\overline v_i=k_i(-\overline \omega_i)$, $\overline u_i=\overline \beta_i$, and choosing the storage function
\[
Z_i:=\frac{\tau_{\alpha, i}}{\rho_i^c} (\alpha_i-\overline \alpha_i)^2+\tau_{\beta, i}(\beta_i-\overline \beta_i)^2\,,
\]
it is easy to see that
\begin{align*}
\dot Z_i &\leq (\frac{1}{\rho^c_i})^2 (v_i-\overline v_i)^2- (\beta_i-\overline \beta_i)^2 \\
&\leq (\frac{\rho^k_i}{\rho^c_i})^2 (\omega_i-\overline \omega_i)^2-  (u_i-\overline u_i)^2\,.
\end{align*}
Therefore the system \eqref{e:example-2nd} has an incremental $L_2$-gain $\leq\frac{\rho^k_i}{\rho^c_i}$, and by Proposition \ref{p:L2}, the equilibrium
$(\overline \varphi, \overline \omega, \overline \xi)$ with $\overline \xi_i=(\overline\alpha_i, \overline\beta_i)$ is asymptotically stable if $\rho_i^k< \rho_i^c D_i$ for each $i\in \calV_g$.
In the special case where $c_i=\frac{1}{2}\alpha_i^2$, $\nabla c_i$ becomes linear, and the stability condition simplifies to $\rho_i^k< D_i$.
The latter is consistent with the result obtained in \cite{zhao2014optimal}.
\EW
\end{example}

\begin{remark}\label{r:diff}
Analogous $L_2$-gain arguments and small gain results were also mentioned in \cite{zhao2014optimal,kasis2017primary}. We have provided the analysis here mainly for two reasons: i) Completeness/concreteness: to provide the explicit form of the Lyapunov functions, and to take into account the subtle technical differences with the model adopted in [13], [14],  such as {the absence} of the damping in the load buses. ii) Comparison: the form of the Lyapunov functions and the $L_2$-gain conditions are provided in order to contrast them with the secant conditions and the corresponding Lyapunov construction in Subsection 4.3.
\end{remark}
\subsection{Incremental secant conditions}\label{ss:secant}
 While $L_2$-gain arguments are powerful and applicable to fairly general classes of nonlinear generation dynamics, they  often lead to conservative conditions that require a substantial amount of damping for stability guarantees. To put forward an alternative approach and obtain less conservative stability conditions, a key observation is that generation dynamics typically can be written as a {\em cascaded} interconnection of output-strictly incrementally passive systems. We impose this observation as an assumption, and will study its applicability on several examples later in the manuscript.

 Before providing the explicit assumption, note that we use the same (incremental) passivity notion as in Definition \ref{d:passive} for systems of ordinary differential equations as a special case. Moreover, we call a static input-output map $y=\phi(u)$ output strictly incrementally passive with respect to a point $(\overline u, \overline y)$, with $\overline y=\phi(\overline u)$, if \eqref{ine:passive} holds with $S=0$ and $Q>0$, i.e.,
\be
\label{ine:passive-static}
0\le
-(y- \overline y)^{\sf T} Q(y-\overline y)+(y-\overline y)^{\sf T} (u-\overline u)\,.
\ee

\begin{assumption}\label{key}
For each $i\in \calV_g$, the system \eqref{e:ui} can be written as a cascaded interconnection of single-input single-output subsystems $\Sigma_{ij}$, $j\in \mathcal{P}=\{1, \ldots, n_\mathcal{P}\}$, such that
\begin{enumerate}
\item Each input-state-output block $\Sigma_{ij}(v_{ij}, \xi_{ij}, z_{ij})$ is output strictly incrementally passive with respect to $(\overline v_{ij}, \overline \xi_{ij}, \overline z_{ij})$, namely \eqref{ine:passive} holds for some storage function $S_{ij}$ and a positive scalar $Q_{ij}$. {The storage function $S_{ij}$ has a strict minimum at
$(\overline v_{ij}, \overline \xi_{ij}, \overline z_{ij})$.}
\item Each static block $\Sigma_{ij}$ given by input-output relation $z_{ij}=\phi_{ij}(v_{ij})$ is
output strictly incrementally passive with respect to $(\overline v_{ij}, \overline z_{ij})$, namely \eqref{ine:passive-static} holds for some positive scalar $Q_{ij}$.
\end{enumerate}
Within the assumption, $\omega_i=v_{i1}$, $z_{i(k-1)}=v_{ik}$ for $k=2,\ldots, n_{\mathcal{P}}$, $z_{in_\mathcal{P}}=u_i$, and $\col(\xi_{ij})$, $j\in \mathcal{P}$, is equal to the vector $\xi_i$ in \eqref{e:ui}, $i\in \calV_g$. The variables with the overlines are defined consistently, noting that $\overline\xi_i=\col(\overline\xi_{ij})$ are such that \eqref{e:feas} is satisfied.
\end{assumption}
\begin{remark}
For linear blocks, the incremental passivity property in Assumption \ref{key} reduces to passivity. For static blocks, the incremental passivity property amounts to an incremental sector boundedness where the slope of nonlinearity does not exceed $Q_{ij}\in \R$.
\end{remark}
Unlike parallel interconnections, cascaded interconnection does not preserve passivity properties, and hence closed loop stability is not automatically guaranteed. 
 However, under Assumption \ref{key}, the shortage of (incremental) passivity can be quantified by adapting the so-called ``secant conditions" \cite{sontag2006passivity,arcak2006diagonal}, to our incremental setting, where loads act as external constant disturbances to the system.  This brings us to the following theorem:
\medskip{}

\begin{theorem}\label{t:secant}
Let Assumptions \ref{a:feas} and \ref{key} hold.
Then, $(\overline \varphi, \overline \omega_g, \overline \xi)$ is an asymptotically stable equilibrium of \eqref{e:varphi-compact} if
\be\label{e:condition}
D_i^{-1}< Q_{i1}\cdots  Q_{in_\mathcal{P}}\, \big(\sec (\frac{\pi}{n_\mathcal{P}+1})\big)^{n_\mathcal{P}+1}\,
\ee
for each $i\in \calV_g$.
\end{theorem}

\medskip{}
\BP
For each $i\in \calV_g$ and $j\in \calP$, let $S_{ij}$ be the storage function obtained from Assumption \ref{key}, where we set $S_{ij}=0$ if $\Sigma_{ij}$ is a static block. Now, for each $i$, let $S_i:=\sum_{j\in \calP} \alpha_{ij}S_{ij}$, where the scalars $\alpha_{ij}\in \R^+$ will be determined afterwards.
In addition, let $z_i=\col(z_{ij})$, $\overline z_i=\col(\overline z_{ij})$, $j\in \mathcal{P}$.
Then, we have
\begin{align}\label{e:dS_i-sec-proof}
\nonumber
\dot S_i&\leq  \sum_{j\in \mathcal{P}} - \alpha_{ij} Q_{ij} (z_{ij}- \overline z_{ij})^2 + \alpha_{ij} (z_{ij}- \overline z_{ij})(v_{ij}-\overline v_{ij}) \\
&= (z_i -\overline z_i)^{\sf T} Z_i (z_i -\overline z_i)+ \alpha_{i1}(z_{i1}-\overline z_{i1})(v_{i1}-\overline v_{i1})\,,
\end{align}
where $Z_i\in \R^{n_\mathcal{P}\times n_\mathcal{P}}$ is a lower triangular matrix with its $(p, q)$th element given by
\[
(Z_i)_{pq}=\begin{cases}
0 & p<q\\
-\alpha_{ip}Q_{ip} & p=q\\
\alpha_{ip} & p=q+1\\
0 & p>q+1\,.
\end{cases}
\]
Here, we have used the fact that $z_{i(k-1)}=v_{ik}$ and $\overline z_{i(k-1)}=\overline v_{ik}$ for $k=2,\ldots, n_{\mathcal{P}}$.
Now, take the Lyapunov function candidate $V:=S^*+\sum_{i\in \calV_g} S_i$ with $S^*$ being equal to the storage function $S$ in \eqref{e:Storage}. Then, by \eqref{e:dS_i-sec-proof} and Proposition \ref{p:passivity}, we obtain that
\begin{align}\label{e:dissip-secant-poof}
\nonumber
&\dot V\leq \\[-2mm]
&\quad\, \sum_{i\in \calV_g} \bbm -\omega_i+\overline \omega_i &\hspace{0mm} (z_i-\overline z_i)^{\sf T}\ebm
\bbm -D_i &  -e_{\mathcal{P}}^{\sf T}   \\[1mm] \alpha_{i1}e_1 & Z_i
\ebm
\bbm -\omega_i+\overline \omega_i \\[1mm] z_i-\overline z_i\ebm,
\end{align}
where $e_1=\bbm 1 & 0 & \cdots & 0\ebm^{\sf T}$, $e_{\mathcal{P}}=\bbm 0 & 0 & \cdots & 1\ebm^{\sf T}$, and we used the fact that $\omega_i=v_{i1}$, $\overline \omega_i=\overline v_{i1}$, $u_i=z_{in_\calP}$, $\overline u_i=\overline z_{in_\calP}$.
Note that
\begin{align*}
\bbm -D_i &  -e_{\mathcal{P}}^{\sf T}   \\[1mm] \alpha_{i1}e_1 & Z_i
\ebm&=
\diag(1, \alpha_{i1}, \dots, \alpha_{in_\mathcal{P}})\\
&\qquad  \times
\bbm
-D_i & 0  & \cdots & 0 & -1\\
1 & -Q_{i1} & \ddots & {} & 0\\
0 & 1 & -Q_{i2} & \ddots & \vdots\\
\vdots & \ddots & \ddots & \ddots & 0\\
0 & \cdots & 0 & 1 & -Q_{in_\mathcal{P}}
\ebm.
\end{align*}
 Then, by \eqref{e:condition} and \cite[Thm. 1]{arcak2006diagonal}, a set of positive scalars  $\alpha_{i1}$, $\ldots$, $\alpha_{in_\mathcal{P}}$ exists such that
 \[
\bbm -D_i &  -e_{\mathcal{P}}^{\sf T}   \\[1mm] \alpha_{i1}e_1 & Z_i
\ebm
+
\bbm -D_i &  -e_{\mathcal{P}}^{\sf T}   \\[1mm] \alpha_{i1}e_1 & Z_i
\ebm^{\sf T} <0\,.
\]
Therefore, by \eqref{e:dissip-secant-poof}, $\dot V$ is nonpositive and is equal to zero whenever $(\omega_i, z_i)=(\overline \omega_i, \overline z_i)$ for all $i\in \calV_g$.
Now, analogous to Proposition \ref{p:L2}, one can invoke the LaSalle's invariance principle and show that the corresponding invariant set of the system, with $\dot{V}=0$, comprises only the equilibrium $(\overline\varphi, \overline\omega, \overline\xi)$.   This completes the proof.
\EP
\medskip{}

\begin{example}\label{ex:deadband}
Suppose that the generation dynamics at each bus $i\in \calV_g$ is given by the first order model
\bse\label{e:example-deadband}
\begin{align}
\tau_{\xi, i}\dot\xi_i&=- \xi_i + k_i(-\omega_i)\\
u_i&=\xi_i
\end{align}
\ese
with $\tau_{\xi, i}\in \R^+$. The map $k_i:\Omega_k \rightarrow \R$ is increasing and satisfies $|k_i(-\omega_i)-k_i(-\overline \omega_i)|\leq  \rho_i |\omega_i-\overline \omega_i|$, $\forall \omega_i, \overline \omega_i\in \Omega_k$, for some $\rho_i\in \R^+$. Typical examples of $k_i$ include deadband nonlinearities and inverse of marginal costs in primary control \cite{zhao2014optimal,kasis2017ifac}. The first order dynamics can be obtained from \cite[Ch. 11]{bergen1999power} by neglecting the fast dynamics of the governor, see e.g. \cite{li2016connecting}. {The dynamics \eqref{e:example-deadband} can also model a decentralized leaky-integral controller \cite{ainsworth2013design,heidari2017ultimate,weitenberg2017robust}, where $k_i(\cdot)$ here is allowed  to be nonlinear.}

Clearly we have
\begin{align}\label{e:k-incpassive}
\nonumber
0&\leq - \frac{1}{\rho_i}|k_i(-\omega_i)-k_i(-\overline \omega_i)|^2 \\
& \qquad - (\omega_i-\overline \omega_i)(k_i(-\omega_i)-k_i(-\overline \omega_i))\,.
\end{align}
Hence, $z_{i1}:=k_i(-\omega_i)$ defines an {output strictly} incrementally passive map. Moreover, by taking the storage function $S_{i2}= \frac{1}{2} \tau_{\xi,i}(\xi_i-\overline \xi_i)^2$ with $\overline \xi_i$ denoting the equilibrium of \eqref{e:example-deadband} resulting from the constant input $-\overline \omega_i$, we have
\[
\dot S_{i2}=-(\xi_i- \overline \xi_i)^2+ (\xi_i- \overline \xi_i)(v_{i2}- \overline v_{i2})\,,
\]
where $v_{i2}=k_i(-\omega_i)=z_{i1}$ and $\overline v_{i2}=k_i(-\overline\omega_i).$ This implies that the system with input-state-output $(v_{i2}, \xi_i, u_i)$ is {output strictly} incrementally passive.  Therefore, Assumption \ref{key} is satisfied with $n_\mathcal{P}=2$, $Q_{i1}=\rho_i^{-1}$ and $Q_{i2}=1$. Consequently, by Theorem \ref{t:secant}, $(\overline \varphi, \overline \omega, \overline \xi)$ is asymptotically stable if
\[
\rho_i <8 D_i\,.
\]
This condition is eight times less conservative than sufficient damping conditions obtained from $L_2$-gain {arguments.} 
\EW
\end{example}

\begin{example}
Let the generation dynamics at each bus $i\in \calV_g$ be given by the nonlinear second-order dynamics in \eqref{e:example-2nd}, see also \cite[Sec. 11.1]{bergen1999power}.
We split the dynamics into three cascaded subdynamics, namely
\begin{align}\label{e:1st-block}
\hspace{-1.32cm}z_{i1}=k_i (-\omega_i)\,,
\end{align}
\vspace{-6mm}
\bse\label{e:2nd-block}
\begin{align}
\tau_{\alpha, i}\dot\alpha_i&=-\nabla c_i(\alpha_i) + v_{i2}\,, \quad v_{i2}=z_{i1}\,,\\
z_{i2}&=\alpha_i\,,
\end{align}
\ese
and
\vspace{-4mm}
\bse\label{e:3rd-block}
\begin{align}
\label{e:3rd-block-beta}
\tau_{\beta, i}\dot\beta_i&= - \beta_i+v_{i3}, \quad v_{i3}=z_{i2}\\
u_i&=\beta_i.
\end{align}
\ese
As before, the first block \eqref{e:1st-block} satisfies the incremental passivity property \eqref{e:k-incpassive} with $\rho_i$ being replaced by $\rho_i^k$.
The storage functions $S_{i2}=\frac{1}{2}\tau_{\alpha_i}(\alpha_i -\overline \alpha_i)^2$ and $S_{i3}=\frac{1}{2}\tau_{\beta,i}(\beta_i -\overline \beta_i)^2$ yields the incremental passivity of the second and third subsystems, \eqref{e:2nd-block} and \eqref{e:3rd-block}, with corresponding coefficients $Q_{i2}=\rho^c_i$ and $Q_{i3}=1$, respectively. Therefore, noting that $n_\mathcal{P}=3$, the secant condition in Theorem \ref{t:secant} reads as
\[
\displaystyle\frac{\rho_i^k}{\rho_c^i}< 4D_i\,.
\]
Again note that the condition above is 4 times less conservative than the one resulting from an $L_2$-gain argument, {see Example \ref{ex:2-order}.}

Next, it is illustrative to consider the same dynamics as before but with an additional nonlinear map at the outputs, namely
\be\label{e:h(u)}
u_i=h_i(\beta_i)\,,
\ee
where $h_i:\Omega_h \rightarrow \R$ is strictly increasing and satisfies  $|h_i(\beta_i)-h_i(\overline \beta_i)|\leq  \rho_i^h |\beta_i-\overline \beta_i|$, $\forall \beta_i, \overline \beta_i\in \Omega_h$, for some $\rho_i^h\in \R^+$.  Hence, $h_i$ defines an incrementally passive map, and can be treated as a new block next to the three subsystems \eqref{e:1st-block}, \eqref{e:2nd-block}, and \eqref{e:3rd-block-beta}. Then, inequality \eqref{e:condition} with $n_{\mathcal P}=4$ gives the stability condition
$\frac{\rho_i^h\rho_i^k}{\rho_c^i}< 2.88 D_i$. However, noting that the secant condition becomes more conservative as the number of cascaded subsystems increases, a compelling alternative is to refine the storage function, and possibly keep the number of cascaded blocks the same. To this end, let $S_{i3}$ be redefined as
\[
S_{i3}:= H_i(\beta_i)-H_i(\overline \beta_i) - \left.(\beta_i- \overline \beta_i) \frac{\partial H_i}{\partial \beta_i}\right|_{\beta_i=\overline \beta_i}\,,
\]
where
\[
H_i(\beta_i)=\tau_{\beta,i}\displaystyle\int_{\overline \beta_i}^{\beta_i} h_i(\tilde \beta_i) d \tilde \beta_i.
\]
We note that $S_{i3}$ is associated with the Bregman distance defined on the function $H_i$ with respect to the point $\overline \beta_i$ \cite{bregman1967relaxation}. Since $h_i$ is strictly increasing, the function $H_i$ is strictly convex, and therefore the storage function $S_{i3}$ is positive definite.  Computing the time derivative of $S_{i3}$ along the solutions of \eqref{e:3rd-block-beta} yields
\begin{align*}
\dot S_{i3}&=-(\beta_i -\overline \beta_i) (h_i(\beta_i)- h_i(\overline \beta_i))\\ &\qquad \quad +(h_i(\beta_i)- h_i(\overline \beta_i))(v_{i3}- \overline v_{i3})\\
&\leq -\frac{1}{\rho_i^h} (u_i- \overline u_i)^2 + (u_i- \overline u_i) (v_{i3}- \overline v_{i3})\,.
\end{align*}
This amounts to the incremental passivity property of \eqref{e:h(u)} with $Q_{i3}=(\rho_i^h)^{-1}$. Hence, Theorem \ref{t:secant} can be applied with $n_{\mathcal P}=3$, which gives the more relaxed stability condition
\[
\frac{\rho_i^h\rho_i^k}{\rho_c^i}< 4 D_i\,. {\tag*{$\square$}}
\]
\end{example}
\begin{remark}\label{r:model-validity}
Note that the proposed results can be used for design purposes as well.
An example is the {``leaky-integral"} controllers \cite{franklin1994feedback,ainsworth2013design,heidari2017ultimate,weitenberg2017robust}, commented in Example \ref{ex:2-order}, where our analysis provides additional flexibility in the design,  and allows to incorporate  turbine-governor dynamics and {practically relevant} nonlinearities such as saturation or deadbands.
It is also worth mentioning that  the proposed analysis can be suitably modified to provide decentralized stability conditions  {when the aforementioned generation dynamics are present in conjunction with some other frequency control schemes that have been proposed in the literature, such as primal-dual algorithms for optimal power sharing} (see e.g. \cite{li2016connecting} and the references therein).  The required modification essentially reduces to adding a (quadratic) term in the proposed Lyapunov functions to compensate for the additional dynamics of the frequency controller.
\end{remark}

\subsection{Exploiting the bounds on line parameters}\label{ss:popov-cond}
Recall that the stability conditions proposed in the previous section are independent of transmission line parameters and are valid for all $\gamma_k=|\beta_{ij}|V_iV_j\in \R^+$ , $k \sim \{i, j\}$, {as long as a synchronous solution exists (see Assumption \ref{a:feas-lin}).}  This feature can be unnecessary if bounds on the transmission line parameters are known. Note that such bounds readily provide bounds on the active power flows due to the boundedness of the sine function.
In this subsection, through a Lyapunov analysis, we investigate conditions under which a synchronous motion of power network {(if exists)} is ``attractive"
for\footnote{{Note that the parameters $\gamma_k$ cannot be arbitrary small, otherwise the active power flow would not be able to compensate for the net-demand at steady-state, see the feasibility condition \eqref{e:feas-linear}.}}
\be\label{e:bound-ki}
\sum_{j\in \calN_i} |\beta_{ij}|V_iV_j  \leq \frac{1}{2}\sigma_i\,,
\ee
given $\sigma_i\in \R^+$, $i\in \calV$.
Toward this end, we make two simplifying assumptions, namely: generation dynamics in \eqref{e:ui} are linear, and we consider aggregated models where each bus has some nonzero inertia meaning that algebraic constraints are absent {(see Remark \ref{r:Hill} on relaxing the latter assumption).}
Note that the overall dynamics are still nonlinear due to the nonlinearity of the power flow.

In this case, the dynamics of $u_i\in \R$ is given by a minimal linear time-invariant system
\begin{align*}
\dot \xi_i&= A_i \xi_i - B_i \omega_i\\
u_i&=C_i \xi_i
\end{align*}
where $A_i\in \R^{n_i\times n_i}$ is invertible, and $B_i\in \R^{n_i\times 1}$, $C_i\in \R^{1\times n_i}$ are nonzero matrices.

\begin{assumption}\label{a:Hurwitz}
The matrix
$$\bbm
-M_i^{-1}D_i & M_i^{-1}C_i\\-B_i & A_i
\ebm$$
 does not have any purely imaginary eigenvalues and $- C_iA_i^{-1}B_i +D_i>0$.
\end{assumption}
\begin{remark}
The assumption on the eigenvalues is required to ensure that the frequency dynamics at each isolated bus, i.e. $\sigma_i=0$, are asymptotically stable.
Note that the condition which will be proposed in Theorem \ref{t:popov}  rules out the possibility of eigenvalues in the open right half plane.
The second condition in Assumption \ref{a:Hurwitz} is a mild assumption on the DC gain of the transfer function from $-\omega_i$ to $u_i$ imposing a negative feedback at steady-state.
\end{remark}
The power network dynamics can be written in vector form as
\bse\label{e:theta-compact-p}
\begin{align}
\dot\theta&=\omega\\
\label{e:swing-popov}
M \dot{\omega}&=-D\omega_i-p(\theta)+p^\ast+u.
\end{align}
\ese
where $p(\theta)=R\Gamma \boldsymbol{\sin}(R^{\sf T}\theta)$ is the vector of power transfer as before. 
The generation dynamics in vector form read as
\bse\label{e:gen-compact-lin}
\begin{align}
\dot \xi&= A \xi - B \omega\\
u&=C \xi
\end{align}
\ese
where the matrices $A$, $B$, $C$, are now block diagonal.

Recall that we are interested in a {\em synchronous motion} of the power network, where the voltage phasors rotate with the same frequency: $\overline \theta_i=\omega^*t+\overline \theta_{0,i}$ for each $i$, with constant $\overline\theta_{0,i}\in \R$. 
Since we are interested in local conditions, under which a synchronous motion is attractive,  the change of coordinates in Section \ref{s:synch} is no longer suitable as it, in general, couples the dynamics of non-adjacent buses.
Therefore, unlike the previous section, here we work with the original coordinates $(\theta, \omega, \xi)$ and  a (time-dependent) synchronous motion $(\overline \theta, \overline \omega, \overline \xi)$. With a little abuse of the notation, we use the set $\calV=\{1, 2, \ldots, n\}$ to denote the set of buses in this subsection.
%
%
%
%

For the model \eqref{e:theta-compact-p}, \eqref{e:gen-compact-lin}, a synchronous motion exists if there exist constant vectors $\overline \theta_0$, $\overline \xi$, and $\overline \omega=\ones\omega^*$ with $\omega^*\in \R$, such that
\bse\label{e:feas-lin}
\begin{align}
0&=-D\ones \omega^* - R\Gamma  \boldsymbol{\sin}(R^{\sf T}\overline\theta_0) + p^*+ C\overline \xi \\
0&=A\overline\xi-B\ones \omega^*
\end{align}
\ese
The condition above can be made more explicit by using the following lemma:
\begin{lemma}\label{l:equib-lin}
The point $(\overline \theta_0, \overline \omega, \overline \xi)$,  with $\overline \omega=\ones \omega^*$, satisfies \eqref{e:feas-lin} if and only if
\[
\overline \xi=A^{-1}B\ones \omega^*, \quad R\Gamma \boldsymbol{\sin}(R^{\sf T}\,\overline \theta_0)= p^*-(D-CA^{-1}B)\ones \omega^*
\]
with
$$
\omega^*=\frac{\ones^{\sf T}p^*}{\ones^{\sf T}(D-CA^{-1}B)^{-1}\ones}.
$$
\end{lemma}
\BP
The proof follows from straightforward algebraic calculations from \eqref{e:feas-lin}.
\EP

By Lemma \ref{l:equib-lin}, existence of a synchronous motion imposes the following feasibility assumption:
\begin{assumption}\textbf{(Existence of a synchronous motion)}\label{a:feas-lin}
There exists a constant vector $\overline \theta_0\in \R^{n}$, with $R^{\sf T}\overline \theta_0 \in (-\frac{\pi}{2}, \frac{\pi}{2})^m,$ such that
\be\label{e:feas-linear}
R\Gamma \boldsymbol{\sin}(R^{\sf T}\,\overline \theta_0)= \big(I_n - \frac{(D-CA^{-1}B)\ones \ones^{\sf T}}{\ones^{\sf T} (D-CA^{-1}B)^{-1}\ones}\big)p^*.
\ee
\end{assumption}

The feasibility of the condition above can be verified using the results available on the solvability of (active) power flow equations, see e.g. \cite{dorfler2013synchronization,jafarpour2017synchronization}. In case the graph $\calG$ is a tree, the incidence matrix $R$ has full column rank and Assumption \ref{a:feas-lin} holds whenever
\[
\norm{\Gamma^{-1}(R^{\sf T}R)^{-1}R^{\sf T} c}_\infty<1,
\]
where $c$ denotes the vector in the right hand side of \eqref{e:feas-linear}.}

The main result of this subsection is stated next, while its proof is postponed to the end of the subsection.

\begin{theorem}\label{t:popov}
Let Assumptions \ref{a:Hurwitz} and \ref{a:feas-lin} hold, and $\sigma_i\in \R^+$ be such that \eqref{e:bound-ki} is satisfied for each $i\in \calV$.
Let $G_i$ denote the transfer function from $p_i^*-p_i(\theta)$ to $\omega_i$, i.e. $G_i(s)=\frac{1}{M_is+C_i(sI-A_i)^{-1}B_i+D_i}$.
Assume that there exists $\rho \in \R^+$, with $-\rho^{-1}$ not being a pole of $G_i$, such that the perturbed transfer matrix
\be\label{e:Hi(s)}
H_i(s):=\sigma_i^{-1}+\frac{1+\rho s}{s}G_i(s)
\ee
is positive real for each $i$.
Then, the vector $(R^{\sf T}\theta, \omega, \xi)$ in \eqref{e:theta-compact-p}, \eqref{e:gen-compact-lin}, locally\footnote{The term locally refers to the fact that solutions are initialized in a suitable neighborhood of the point
$(R^{\sf T}\overline\theta, \overline \omega, \overline \xi)$.}converges to $(R^{\sf T}\overline\theta, \overline \omega, \overline \xi)$. Such convergence is established by the Lyapunov function $W+Z$ with $W$ and $Z$ given by \eqref{e:W} and \eqref{e:Z}, respectively.
\end{theorem}

The result is inspired by the classical Popov criterion, with three notable differences: i) An immediate application of the Popov criterion on the networked dynamics results in fully centralized conditions, whereas the conditions here are primarily local (see Remark \ref{r:local}). ii)  The Popov criterion is stated in terms of {\em strict} positive realness of a perturbed transfer function \cite[Ch. 7]{khalil2002nonlinear}, while the result here is provided in terms of positive realness only. This allows us to cope with the presence of the pure integrator in $H_i$, which would otherwise be difficult to remove with a local perturbation argument. The challenge imposed by the lack of {\em strict} passivity in $H_i$ will be overcome by studying asymptotic behavior of the system using Barbalat's Lemma. iii) Due to the presence of the term $p^*$, acting as a constant disturbance to \eqref{e:theta-compact-p}, suitable incremental Lyapunov functions are needed to establish convergence of the solutions to a synchronous motion, see also Remark \ref{r:Lyapunov}.

\begin{remark}\label{r:lmi}
The positive realness condition in Theorem \ref{t:popov} can be equivalently expressed in the state-space domain using matrix inequalities \cite{anderson1967system,willems1972dissipative}.  The additional technical assumption $-\rho^{-1}$ not being a pole of $G_i$ is then translated to $-\rho^{-1}$ not being an eigenvalue of the matrix in Assumption \ref{a:Hurwitz}. The latter is necessary to ensure the existence of a {\em positive definite} solution to the aforementioned matrix inequalities, see also Lemma \ref{l:minimality}.
\end{remark}

\begin{remark}\label{r:local}
Note that, with the exception of the constant $\rho$, only local/distributed information is exploited in the condition of Theorem \ref{t:popov}. More precisely, we rely on three sorts of information: i) Nodal information, which involves knowing the transfer matrix $G_i$, or in other words the matrices $M_i$, $D_i$, $A_i$, $B_i$, and $C_i$. ii) Neighboring information, associated with the bound $\sigma_i$ in \eqref{e:bound-ki}. iii) Global information, which accounts for the parameter $\rho$.  The latter dictates a protocol that must be followed by each bus such that stability is not jeopardized by the interconnection via the power transfers. If the network is expanded,
stability can be guaranteed providing that the newly added buses satisfy the same protocol.
In the special case where $G_i$ is passive, for each $i$, the condition in Theorem \ref{t:popov} becomes independent of $\rho$, by taking the limit of $H_i$ as $\rho$ tends to infinity. In that case, the constant scalars $\sigma_i$ in \eqref{e:bound-ki} can be chosen arbitrary large as expected.
\end{remark}

\begin{remark}\label{r:rho}
	The positive realness condition in Theorem \ref{t:popov} holds if $G_i$ has no poles on the closed right half plane, and for each $i$ we have that
	\begin{equation}\label{e:freq-cond}
	\sigma_i^{-1}+\rho X_i(\nu)+\frac{Y_i(\nu)}{\nu}>0, \qquad \forall \nu >0,
	\end{equation}
	where $X_i(\nu)=\Re(G_i(j\nu))$ and $Y_i(j\nu)=\Im(G_i(j\nu))$, $X_i(0)>0$.
	Note that the ratio $\frac{Y_i(\nu)}{\nu}$ is bounded and  converges to zero as $\nu$ tends to infinity.
	At the low frequencies, we have $X_i(\nu)>0$ and  hence there exists $\underline{\rho_i}\geq 0$ such that \eqref{e:freq-cond} is satisfied for all $\rho > \underline{\rho_i}$.
	On the other hand, at higher frequencies where $X_i(\nu)$ may no longer be positive, but the ratio $\frac{Y_i(\nu)}{\nu}$ becomes small,  one should choose
	$\rho<\overline \rho_i$ for some appropriately chosen  $\overline \rho_i>0$. Consequently, in order to satisfy \eqref{e:freq-cond}, it must hold that $\underline{\rho_i} < \overline{\rho_i}$, and the intervals  $(\underline{\rho_i}, \overline \rho_i)$, $i\in \mathcal{I}$, should have a nonempty intersection. Note that if $G_i$ is passive, then $\overline{\rho_i}$ can be chosen arbitrary large. Loosely speaking, the existence and the corresponding value of $\rho$ satisfying  \eqref{e:freq-cond} will be determined by the buses dynamics that are farthest away from passivity and are strongly coupled to the rest of the network.
\end{remark}

\begin{remark}\label{r:Lyapunov-incrementsl}
Note that the vector $p^*$ contains information on the loads, which may not be accurately available. The incremental construction of Lyapunov functions pursued here gives rise to stability certificates that are independent of $p^*$, as long as the model \eqref{e:theta-compact-p} is valid. Notice that the vector $p^*$ only contributes to the feasibility condition \eqref{e:feas-linear}, and does not appear in \eqref{e:Hi(s)}.
\end{remark}

\begin{remark}\label{r:Lyapunov}
The nonlinear Lyapunov analysis carried out here provides in general a larger region of attraction, compared to the one obtained from linearization. Interestingly, it can be verified that substituting $\sin(\delta)$ and $\cos(\delta)$, $\delta\in \R$, in \eqref{e:Z}, by their second degree Taylor polynomials, namely $\delta$ and $1-\frac{\delta^2}{2}$, respectively, yields a quadratic Lyapunov function which can be used to establish stability properties of the linearized model. The non-quadratic Lyapunov function exploited here, or in other words the full Taylor series of sine and cosine functions, provides additional flexibility that are used to cope with the nonlinearity of the power flows.
It should be noted though that linearizing the power flow could facilitate the use of input/output approaches \cite{devane2017primary}, which can allow the parameter $\rho$ in \eqref{e:Hi(s)} to be an expression in the frequency-domain.
An investigation of the underlying structure of the Lyapunov functions in such cases, and the nonlinearities they could efficiently capture  is an interesting problem and a part of ongoing work.
\end{remark}

	\begin{remark}\label{r:Hill}
		The result of Theorem \ref{t:popov} can be extended to a structure-preserving model, where the load buses are given by \cite{bergen1981structure}
		\[
		D_i\dot\theta_i=-p_i(\theta)+p_i^*, \quad i\in \mathcal{V}_L.
		\]
		In this case, the positive realness condition in Theorem \ref{t:popov} needs to be verified only for the generation buses. The Lyapunov function that establishes the stability result for the structure-preserving case is given by
		\[
		\hat W+ Z + \frac{1}{2}\sum_{i \in \mathcal{V}_\ell} D_i(\theta_i- \overline \theta_i)^2,
		\]
		where $\hat W$ has the same expression as $W$ in \eqref{e:W} but with $i\in \calV_g$, and the function $Z$ is given by \eqref{e:Z}.
	\end{remark}

\begin{example}\label{ex:popov}
Consider a four area power network whose dynamics are governed by \eqref{e:theta-compact-p}, see \cite{nabavi2013topology} on how a four area network equivalent can be obtained for the IEEE New England 39-bus system or the South Eastern Australian 59-bus. Suppose that the generation dynamics are given by the second-order system
\begin{align*}
\tau_{\alpha, i}\dot\alpha_i&=-\alpha_i - k_i\omega_i\\
\tau_{\beta, i}\dot\beta_i&= - \beta_i+\alpha_i\\
u_i&=\beta_i\,,
\end{align*}
and we set $\tau_{\alpha,i}=0.5$ and $\tau_{\beta,i}=1$ for each $i$, and the voltage magnitudes are $V_i\simeq 1$(pu). The inertia, damping, and droop gains of the areas are provided in Table \ref{t:param}. Note that, for illustration purposes, the numerical value of the droop gain $k_1>0$ has not been fixed.  To evaluate the condition in Theorem \ref{t:popov}, {the remaining required parameters are} the values of $\sigma_i$, $i=1, 2, 3, 4.$ Suppose that $\sigma_1\geq \max (\sigma_2, \sigma_3, \sigma_4)$. Then, the proposed stability condition can be verified given the pair $(k_1, \sigma_1)$. For different values of $\sigma_1$(pu), the maximum droop gain $k_1$ for which the stability certificates of Theorem \ref{t:popov} hold are provided in Table \ref{tab:popov}.

\begin{table}[t]
	\centering
	\caption{Simulation parameters}
	\label{t:param}
	\scalebox{1}{
	\begin{tabular}{cccccccccccc}
		\toprule
		{\rm Areas} & 1 & 2 & 3 & 4\\[1.2mm]
		\midrule
		$M_i$ & 5.5 & 3.98 & 4.49 & 4.22\\[1.2mm]
		$D_i$ & 1.60 & 1.22 & 1.38 & 1.42\\[1.2mm]
		$k_{i}$ &  $k_1$ & 7 & 8 & 9\\[1.2mm]
		\bottomrule
	\end{tabular}}
	\vspace{0.2cm}
\end{table}

\begin{table}[h]
	\centering
	\caption{Stability certificates in Example \ref{ex:popov}}
	\label{tab:popov}
	\scalebox{1}{
	\begin{tabular}{ccccccc}
		\toprule
		$\sigma_1$ & 0 & 5 & 10 & 15 & 20  & $\geq\,$30\\[1.1mm]
		$k_{1}$ & 24.3  & 20 & 16.9 & 15.2 & 14.3 & $\simeq$\,13.9 \\[1.1mm]
		\bottomrule
	\end{tabular}}
	\vspace{0.1cm}
\end{table}

As can be seen from the table, at $\sigma_1=0$, which corresponds to the isolated bus dynamics, the maximum allowed droop gain is
$k_1=24.3$. As the strength of the coupling, i.e. $\sigma_1$, increases, the value of $k_1$ that can be tolerated in view of stability decreases. Eventually, for $\sigma\geq 30$, we have $k_1\simeq 13.9$. In fact, the latter corresponds to the special case where bus dynamics are passive. It is worth mentioning that an application of secant conditions returns $k_1\leq 8D_1=12.8$ in this case.
This is expected as the secant conditions allow for nonlinear droop gains, and are independent of the inertia, time constants of the model, {and most importantly the bounds on the line parameters}.
\EW
\end{example}

The rest of this subsection is dedicated to the proof of the main result.

\noindent\textit{Proof of Theorem \ref{t:popov}:} Let $z_i:=\sigma_i^{-1}(p_i^*-p_i)+\theta_i +\rho \omega_i$ for each $i\in \calV$.  Then, clearly $H_i$ is the transfer function from $p_i^*-p_i$ to $z_i$.
The transfer function $H_i$ admits the state space realization
\bse\label{e:Hi}
\begin{equation}
\bbm
\dot\theta_i\\
\dot\omega_i\\
\dot\xi_i
\ebm
=
\bbm
0&1&0\\0& -M_i^{-1}D_i & M_i^{-1}C_i\\0 & -B_i & A_i
\ebm
\bbm
\theta_i\\
\omega_i\\
\xi_i
\ebm+
\bbm
0\\
M_i^{-1}\\
0\\
\ebm
v_i,
\end{equation}
\begin{equation}
z_i= \bbm 1\,& \rho \,& 0 \,\ebm \bbm \theta_i \\ \omega_i \\ 0 \ebm+ \sigma_i^{-1} v_i\,,
\end{equation}
\ese
where $v_i:=p_i^*-p_i(\theta)$.
More compactly, we denote the realization above as
\[
\dot x_i=\calA_i x_i + \calB_i v_i\,, \quad z_i=\calC_i x_i + \sigma_i^{-1} v_i\,,
\]
where $x_i=\col(\theta_i, \omega_i, \xi_i)$.
The realization above is minimal as shown in the following lemma. The proof is straightforward, yet is provided in Appendix for the sake of completeness.
\begin{lemma}\label{l:minimality}
Let Assumption \ref{a:Hurwitz} hold, and assume that $-\rho^{-1}$ is not a pole of $G_i$. Then, the pair $(\calA_i, \calB_i)$ is controllable and the pair $(\calC_i, \calA_i)$  is observable.
\end{lemma}
%
%
%
\noindent\textit{Proof of Theorem \ref{t:popov} (continued):} Since $H_i$ is positive real and \eqref{e:Hi} is minimal, there exists a quadratic storage function $W_i(x_i)=x_i^{\sf T}X_ix_i$ with $X_i>0$ such that
$\dot{W}_i\leq z_i^{\sf T}v_i$. By linearity,
$W_i(x_i-\overline x_i)=(x_i-\overline x_i)^{\sf T}X_i(x_i- \overline x_i)$ satisfies $\dot{W}_i\leq
(z_i-\overline z_i)^{\sf T}(v_i-\overline v_i)$, where $\overline x_i:=(\overline \theta_i, \overline \omega_i, \overline \xi_i)$, and $\overline z_i=\calC_i \overline x_i+\sigma_i^{-1}\overline v_i$. This amounts to the incremental passivity property of \eqref{e:Hi}. Let
\be\label{e:W}
W(x-\overline x):=\sum_{i\in \calV} W_i(x_i-\overline x_i)=\sum_{i\in \calV} (x_i-\overline x_i)^{\sf T}X_i(x_i- \overline x_i)\,.
\ee
Then, in vector form, we have
\begin{align}\label{e:Wdot}
\nonumber
\dot{W}&\leq (v-\overline v)^{\sf T}(z-\overline z)\\
\nonumber
&=-\big(p(\theta)- p(\overline \theta)\big)^{\sf T} \\
\nonumber
&\qquad \qquad \big(\theta-\overline \theta + \rho(\omega - \overline \omega) - \Sigma^{-1}
(p(\theta)-p(\overline \theta))\big)\\
\nonumber
&=-\big(\Gamma \boldsymbol{\sin}(\eta)-\Gamma \boldsymbol{\sin}(\overline \eta)\big)^{\sf T}\\
\nonumber
&\qquad \big((\eta-\overline \eta)- R^{\sf T}\Sigma^{-1} R (\Gamma \boldsymbol{\sin}(\eta)-\Gamma \boldsymbol{\sin}(\overline \eta))\big)\\
&\quad \,-\rho \big(\Gamma \boldsymbol{\sin}(\eta)-\Gamma \boldsymbol{\sin}(\overline \eta)\big)^{\sf T} R^{\sf T}(\omega-\overline \omega)\,,
\end{align}
where $\eta:=R^{\sf T}\theta$ and $\overline \eta:=R^{\sf T}\overline \theta=R^{\sf T}\theta_0$, and $\Sigma:=\diag(\sigma_i)$.
To proceed further, we need the following algebraic result, whose proof is provided in Appendix.

\begin{lemma}\label{l:algeb}
It holds that $\Gamma^{-1}-R^{\sf T}\Sigma^{-1}R\geq 0$.
\end{lemma}

\noindent\textit{Proof of Theorem \ref{t:popov} (continued):}

By Lemma \ref{l:algeb} and \eqref{e:Wdot}, we obtain that
\begin{align*}
\dot W&\leq -\big(\Gamma \boldsymbol{\sin}(\eta)-\Gamma \boldsymbol{\sin}(\overline \eta)\big)^{\sf T}\\
&\qquad \qquad \big((\eta-\overline \eta)+\rho R^{\sf T}\omega-  (\boldsymbol{\sin}(\eta)- \boldsymbol{\sin}(\overline \eta))\big)\,,
\end{align*}
where we also used the fact that $R^{\sf T}\overline \omega=0$.
Now, we define the Bregman distance type function  \cite{bregman1967relaxation}
\begin{align}\label{e:Z}
\nonumber
Z(\theta, \overline \theta)&:= \rho\big(-\ones^{\sf T} \Gamma{\boldsymbol{\cos}}(R^{\sf T}\theta)+ \ones^{\sf T} \Gamma {\boldsymbol{\cos}}(R^{\sf T}\overline \theta)\big)\\
&\qquad -\rho (\theta -\overline \theta)^{\sf T}R\Gamma \boldsymbol{\sin}(R^{\sf T}\overline \theta)\,,
\end{align}
where $\boldsymbol{\cos}(\cdot)$ is interpreted element-wise. The function above is nonnegative for $R^{\sf T}\theta\in (-\frac{\pi}{2}, \frac{\pi}{2})^m$ and is equal to zero whenever $R^{\sf T}\theta=R^{\sf T}\overline \theta$. Note that $Z$ does not  explicitly depend on time (other than via the system states) as
$R^{\sf T}\overline \theta=R^{\sf T}\overline \theta_0$ is constant. Computing the time derivative of $Z$ along the solutions of the system yields
\[
\dot Z=\rho(R\Gamma \boldsymbol{\sin}(\eta)- R\Gamma \boldsymbol{\sin} (\overline \eta))^{\sf T} \omega\,,
\]
where again $\eta=R^{\sf T}\theta$ and $\overline \eta=R^{\sf T}\overline \theta$. Therefore, by defining $V:=W+Z$, we obtain that
\begin{align}\label{e:dotS-exp}
\nonumber
\dot V&=\dot W+ \dot Z\\
\nonumber
&\leq -\big(\boldsymbol{\sin}(\eta)- \boldsymbol{\sin}(\overline \eta)\big)^{\sf T}\Gamma \big(\eta-\overline \eta - (\boldsymbol{\sin}(\eta)- \boldsymbol{\sin}(\overline \eta)) \big)\\
\nonumber
&=-\sum_{k\sim \{i, j\}}  \gamma_k\big({\sin}(\eta_k)- {\sin}(\overline \eta_k)\big)\\[-3mm]
&\qquad \qquad \qquad \quad \big(\eta_k-\overline \eta_k - ({\sin}(\eta_k)-{\sin}(\overline \eta_k)) \big)\,.
\end{align}
By the mean value theorem, we have
\begin{align}\label{e:eta-mean}
\nonumber
&\gamma_k\big({\sin}(\eta_k)- {\sin}(\overline \eta_k)\big)\big(\eta_k-\overline \eta_k - ({\sin}(\eta_k)-{\sin}(\overline \eta_k))\\
&\quad =\gamma_k(\eta_k-\overline \eta_k)^2 \cos(\tilde{\eta}_k) (1-\cos(\tilde{\eta}_k))\,,
\end{align}
for some $\tilde{\eta}_k$ which can be written as a convex combination of $\eta_k$ and $\overline \eta_k$. Therefore, the time derivative of $V$ is nonpositive whenever $\eta=R^{\sf T}\theta \in (-\frac{\pi}{2}, \frac{\pi}{2})^m$.

Now, suppose that the vector $\eta=R^{\sf T}\theta(t)$ belongs to a closed subset of $(-\frac{\pi}{2}, \frac{\pi}{2})^m$ for all time.
{This is always possible by initializing $R^T\theta$ sufficiently close to the point $R^{\sf T}\overline \theta=R^{\sf T}\overline \theta_0$, noting that $Z$ is positive definite with respect to $\eta=R^{\sf T}\theta$, and that $R^{\sf T}\overline \theta \in (-\frac{\pi}{2}, \frac{\pi}{2})^m.$}
Notice that $V$ explicitly depends on time due to the term containing $\overline \theta$ in $W$, however the right hand side of \eqref{e:dotS-exp} is only a function of states bearing in mind that the vector $\overline \eta=R^{\sf T}\overline \theta$ is constant.
By integrating both sides of \eqref{e:dotS-exp}, and noting that $V$ is nonnegative, we have
$$
\int _0^\infty \phi(\tau) d\tau \leq V(x(0), \overline x(0))\,,
$$
where $-\phi(\eta)$ denotes the right hand side of \eqref{e:dotS-exp}. Noting that $\phi(\eta)$ is nonnegative, the integral on the left hand side of the inequality above is well-defined. Baring in mind that $Z$ does not explicitly depend on time and is positive definite with respect to $\eta=R^{\sf T}\theta$, the Lyapunov function $V$ is bounded from below by a positive definite function of  $(\eta, \omega, \xi)$ which does not explicitly depend on time. Therefore, recalling that $\dot V$ is nonpositive, we
have that $\eta$, $\omega$, and $\xi$ are bounded.
Then, the time derivative of $\phi(\eta)$ is bounded, and thus $\phi$ is uniformly continuous. By exploiting Barbalat's Lemma \cite[Lem.8.2.]{khalil2002nonlinear}, we then obtain that $\lim_{t\rightarrow \infty} \phi(\eta(t))=0$. As $\eta$ belongs to a closed subset of $(-\frac{\pi}{2}, \frac{\pi}{2})^m$, by \eqref{e:eta-mean} we find that $\lim_{t\rightarrow \infty} \eta(t)=\overline \eta.$ By Assumption \ref{a:Hurwitz} and positive realness of $H_i$, the dynamics from $v_i-\overline v_i$ to $(\omega_i-\overline \omega_i, \xi- \overline \xi_i)$ are given by a linear asymptotically stable system, see \eqref{e:Gi}. Consequently, as $\eta$ converges to $\overline \eta$ and $v$ to $\overline v$, we conclude that $\lim_{t\rightarrow \infty} \omega(t)=\overline \omega$ and $\lim_{t\rightarrow \infty} \xi(t)=\overline \xi$. This completes the proof.
\EP
%

\section{Conclusions}\label{s:conclusions}
We have 
provided a {Lyapunov stability} analysis of a differential algebraic model of frequency dynamics in power networks with turbine governor dynamics, static and dynamic nonlinearities.
In particular, we have shown that 
secant gain conditions which rely on suitable cascaded decomposition of the generation dynamics, can lead to decentralized stability conditions with reduced conservatism.
Furthermore, for linear generation dynamics, we have derived Popov-like conditions
that reduce the conservatism even further, by exploiting additional local information associated with the coupling strength among the bus dynamics.
Numerical examples illustrate that the latter conditions provide improvements in the case the bus dynamics are weakly coupled. As expected, these also coincide with conditions obtained from passivating the bus dynamics as the coupling strength tends to infinity. Interesting directions for future research are to include voltage control dynamics, secondary frequency control schemes, extensions to lossy networks, as well as the use of more involved classes of Lyapunov functions that can provide further flexibility in the analysis.

\section*{Appendix}

\begin{pfof}{Lemma \ref{l:minimality}}
By PBH controllability test, the pair  $(\calA_i, \calB_i)$ is controllable if and only if the matrix
\be\label{e:ctrb}
\bbm -\lambda & 1 & 0\\0 & -B_i & A_i-\lambda I\ebm
\ee
is full row rank for all $\lambda \in \C$. Suppose that there exists a vector $\zeta=\col(\zeta_1, \zeta_2)$ belongs to the left kernel of the matrix in \eqref{e:ctrb}.
We distinguish between the two cases $\lambda =0$ and $\lambda\neq 0$.
First, let $\lambda =0$. Then, we have $\zeta_2^{\sf T}A_i=0$, which implies that $\zeta_2=0$ noting that $A_i$ is nonsingular. This results in $\zeta_1=0$, and thus in controllability of $(\calA_i, \calB_i)$. Now, consider the case where $\lambda \neq 0$. Then, $\zeta_1$ is necessarily zero, and we obtain $\zeta_2^{\sf T} \bbm -B_i & A_i-\lambda I\ebm=0$.
By controllability of $(A_i, -B_i)$, we conclude that $\zeta_2=0$, and hence $(\calA_i, \calB_i)$ is controllable.

For observability of $(\calC_i, \calA_i)$, we need to show that the matrix $\mathcal{Q}_i=\bbm \calA_i^{\sf T}-\lambda I & \calC_i^{\sf T}\ebm^{\sf T}$ has full column rank for all $\lambda\in \C$.
Suppose that there exists a vector $\zeta= \col(\zeta_1, \zeta_2, \zeta_3)$ in the (right) kernel of $\mathcal{Q}_i$, where the partitioning is in accordance with \eqref{e:Hi}.  Then, we have $\zeta_2=\lambda \zeta_1$ and $(1+\lambda \rho)\zeta_1=0$. If $\lambda=0$, then we obtain that $\zeta_1=0$, $\zeta_2=0$, and the result  follows from observability of $(C_i, A_i)$.
If both $\lambda$ and $(1+\rho \lambda)$ are nonzero, then we find again that $\zeta_1=0$ and $\zeta_2=0$ which results in observability of $(C_i, A_i)$.

Finally, note that
\bse\label{e:Gi}
\begin{align}
\bbm
\dot\omega_i\\
\dot\xi_i
\ebm
&=
\bbm
-M_i^{-1}D_i & M_i^{-1}C_i\\-B_i & A_i
\ebm
\bbm
\omega_i\\
\xi_i
\ebm+
\bbm
M_i^{-1}\\
0\\
\ebm
v_i \\
y_i&= \omega_i\,,
\end{align}
\ese
gives a minimal realization of $G_i$. In fact, it is easy to see that the controllability property follows from controllability of the pair $(A_i,- B_i)$, and the observability property is deduced from observability of $(C_i, A_i)$. Hence, the fact that $-\rho^{-1}$ is not a pole of $G_i$ implies that $-\rho^{-1}$ is not an eigenvalue of $\calA_i^\omega$, where $\calA_i^\omega$ denotes the state matrix in \eqref{e:Gi}.
Clearly, it suffices to check the rank of $\mathcal{Q}_i$ for all $\lambda\in \sigma(\calA_i)=\sigma(\calA^\omega_i)\cup \{0\},$
where $\sigma(\cdot)$ denotes the spectrum of the matrix, and we used the block triangular structure of $\calA_i$ to write the last equality. Hence, the fact that $-\rho^{-1}\notin \sigma(\calA^\omega_i)$ yields
$1+\rho\lambda \neq0$, for all $\lambda\in \sigma(\calA_i)$. This completes the proof, since observability under the condition $1+\rho\lambda \neq0$ was established before.
\end{pfof}

\begin{pfof}{Lemma \ref{l:algeb}}
Let the matrix $L$ be defined as $L:=R\Gamma R^{\sf T}$. Note that $L$ is a Laplacian matrix by construction.
First, we show that the eigenvalues of the matrix $\Sigma^{-1}L$ are not greater than $1$. As this matrix is similar to $\Sigma^{-\frac{1}{2}}L\Sigma^{-\frac{1}{2}}$, its eigenvalues are real and nonnegative. By Gershgorin circle theorem, it is easy to see that the eigenvalues of $\Sigma^{-1}L$ are not greater than $1$ if $2\sigma_i^{-1}L_{ii}\leq1$ for each $i$.
The latter inequality holds since, by \eqref{e:bound-ki}, $\sigma_i\geq 2\sum_{j\in \calN_i} \beta_{ij}V_iV_j=2L_{ii}$.
Now, noting that $\Sigma^{-1}L$ is similar to $\Sigma^{-\frac{1}{2}}R\Gamma R^{\sf T}\Sigma^{-\frac{1}{2}}$, it shares the same nonzero eigenvalues as the matrix $\Gamma^{\frac{1}{2}} R^{\sf T}\Sigma^{-1}R \Gamma^{\frac{1}{2}}$. 
Therefore, the eigenvalues of the latter matrix are not greater than $1$ either, and we have
\[
I-\Gamma^{\frac{1}{2}} R^{\sf T}\Sigma^{-1}R \Gamma^{\frac{1}{2}}\geq 0\,.
\]
Finally, by a congruent transformation, the inequality above is equivalent to $\Gamma^{-1}-R^{\sf T}\Sigma^{-1}R\geq 0$.
\end{pfof}

\bibliographystyle{IEEEtran}
\bibliography{nima2}

\end{document}